\documentclass{amsart}
\usepackage{amssymb}

\usepackage{amsmath}
\usepackage{amscd}
\usepackage{thmdefs}

\input{tcilatex}
\theoremstyle{definition}
\theoremstyle{remark}
\numberwithin{equation}{section}

\input tcilatex

\begin{document}
\title[Projections in Graph $W^{*}$-Algebras]{Projections in Graph $W^{*}$-Algebras}
\author{Ilwoo Cho}
\address{Univ. of Iowa, Dep. of Math, Iowa City, IA, U. S. A}
\email{ilcho@math.uiowa.edu}
\date{}
\dedicatory{}
\thanks{}
\maketitle

\begin{abstract}
In this paper, we will consider the projections in a graph $W^{*}$-algebra $%
W^{*}(G).$ Let $T=L_{w_{1}}^{u_{w_{1}}}...L_{w_{n}}^{u_{w_{n}}}$ be a fixed
operator in $W^{*}(G),$ where $w_{1},$ $...,$ $w_{n}$ $\in $ $\Bbb{F}^{+}(G)$
and $u_{w_{1}},$ $...,$ $u_{w_{n}}$ $\in $ $\{1,$ $*\}.$ We will find the
conditions when this operator $T$ is a projection. In other words, we will
characterize the cases when $T=L_{v},$ for some $v\in V(G).$ In this paper,
we show that there exists a vertex $v$ such that $T=L_{v}$ if and only if
the lattice path of $T$ has the $*$-axis-property if and only if $%
w_{1}^{t_{1}}$ $...$ $w_{n}^{t_{n}}$ $=$ $v$ in the free groupoid $\Bbb{F}%
(G),$ where $t_{j}=1$ if $u_{j}=1$ and $t_{j}=-1$ if $u_{j}=*,$ for all $%
j=1,...,n.$
\end{abstract}

\strut

In [2], Kribs and Power defined the free semigroupoid algebras and obtained
some properties of them. Roughly speaking, graph $W^{*}$-algebras are $W^{*}$%
-topology closed version of free semigroupoid algebras. Throughout this
paper, let $G$ be a countable directed graph and let $\mathbb{F}^{+}(G)$ be
the free semigroupoid of $G,$ in the sense of Kribs and Power. i.e., it is a
collection of all vertices of the graph $G$ as units and all admissible
finite paths, under the admissibility. As a set, the free semigroupoid $%
\mathbb{F}^{+}(G)$ can be decomposed by

\strut

\begin{center}
$\mathbb{F}^{+}(G)=V(G)\cup FP(G),$
\end{center}

\strut

where $V(G)$ is the vertex set of the graph $G$ and $FP(G)$ is the set of
all admissible finite paths. Trivially the edge set $E(G)$ of the graph $G$
is properly contained in $FP(G),$ since all edges of the graph can be
regarded as finite paths with their length $1.$ We define a graph $W^{*}$%
-algebra of $G$ by

\strut

\begin{center}
$W^{*}(G)\overset{def}{=}\overline{%
\mathbb{C}[\{L_{w},L_{w}^{*}:w\in
\mathbb{F}^{+}(G)\}]}^{w},$
\end{center}

\strut

where $L_{w}$ and $L_{w}^{*}$ are creation operators and annihilation
operators of $\xi _{w}$ on the generalized Fock space $H_{G}=l^{2}\left( %
\mathbb{F}^{+}(G)\right) $ induced by the given graph $G,$ respectively (See
[1]). Let $w_{1},...,w_{n}\in \Bbb{F}^{+}(G)$ and $u_{1},...,u_{n}\in
\{1,*\}.$ Define

\strut

\begin{center}
$T=L_{w_{1}}^{u_{1}}...L_{w_{n}}^{u_{n}}\in W^{*}(G).$
\end{center}

\strut

In this paper, we want to find the condition when $T$ is a projection. i.e.,
we want to find the condition when there exists $v\in V(G)$ such that $%
T=L_{v}.$ In [1], we already showed that $T=L_{v}$ if and only if the
lattice path of $T$ has the so-called the $*$-axis-property with $*=v.$ In
this paper, we will find one more condition which is equivalent to the $*$%
-axis-property of $T.$ We can show that $T=L_{v}$ is a projection if and
only if the corresponding element $w_{1}^{t_{1}}$ $...$ $w_{n}^{t_{n}}$ $\ $%
is nothing but $v$ in the free groupoid $\Bbb{F}(G),$ induced by the
shadowed graph $G^{\symbol{94}}=G\cup G^{-1}$ of the given graph $G,$ with a
certain relation on its free semigroupoid $\Bbb{F}^{+}(G^{\symbol{94}}).$
Here, the shadowed graph $G^{\symbol{94}}$ is the countable directed graph
with

$\strut $

\begin{center}
$V(G^{\symbol{94}})=V(G)$ and $E(G^{\symbol{94}})=E(G)\cup E(G^{-1})$,
\end{center}

\strut

where $G^{-1}$ is the shadow of the graph $G.$

\strut

\strut

\strut

\section{Graph $W^{*}$-Algebras}

\strut

\strut

Let $G$ be a countable directed graph and let $\Bbb{F}^{+}(G)$ be the free
semigroupoid of $G.$ i.e., the set $\mathbb{F}^{+}(G)$ is the collection of
all vertices as units and all admissible finite paths of $G.$ Let $w$ be a
finite path with its source $s(w)=x$ and its range $r(w)=y,$ where $x,y\in
V(G).$ Then sometimes we will denote $w$ by $w=xwy$ to express the source
and the range of $w.$ We can define the graph Hilbert space $H_{G}$ by the
Hilbert space $l^{2}\left( \mathbb{F}^{+}(G)\right) $ generated by the
elements in the free semigroupoid $\mathbb{F}^{+}(G).$ i.e., this Hilbert
space has its Hilbert basis $\mathcal{B}=\{\xi _{w}:w\in \mathbb{F}%
^{+}(G)\}. $ Suppose that $w=e_{1}...e_{k}\in FP(G)$ is a finite path with $%
e_{1},...,e_{k}\in E(G).$ Then we can regard $\xi _{w}$ as $\xi
_{e_{1}}\otimes ...\otimes \xi _{e_{k}}.$ So, in [2], Kribs and Power called
this graph Hilbert space the generalized Fock space. Throughout this paper,
we will call $H_{G}$ the graph Hilbert space to emphasize that this Hilbert
space is induced by the graph.

\strut

Define the creation operator $L_{w},$ for $w\in \mathbb{F}^{+}(G),$ by the
multiplication operator by $\xi _{w}$ on $H_{G}.$ Then the creation operator 
$L$ on $H_{G}$ satisfies that

\strut

(i) \ $L_{w}=L_{xwy}=L_{x}L_{w}L_{y},$ for $w=xwy$ with $x,y\in V(G).$

\strut

(ii) $L_{w_{1}}L_{w_{2}}=\left\{ 
\begin{array}{lll}
L_{w_{1}w_{2}} &  & \text{if }w_{1}w_{2}\in \mathbb{F}^{+}(G) \\ 
&  &  \\ 
0 &  & \text{if }w_{1}w_{2}\notin \mathbb{F}^{+}(G),
\end{array}
\right. $

\strut

\ \ \ \ for all $w_{1},w_{2}\in \mathbb{F}^{+}(G).$

\strut

Now, define the annihilation operator $L_{w}^{*},$ for $w\in \mathbb{F}%
^{+}(G)$ by

\strut

\begin{center}
$L_{w}^{\ast }\xi _{w^{\prime }}\overset{def}{=}\left\{ 
\begin{array}{lll}
\xi _{h} &  & \text{if }w^{\prime }=wh\in \mathbb{F}^{+}(G)\xi \\ 
&  &  \\ 
0 &  & \text{otherwise.}
\end{array}
\right. $
\end{center}

\strut

The above definition is gotten by the following observation ;

\strut

\begin{center}
$
\begin{array}{ll}
<L_{w}\xi _{h},\xi _{wh}>\, & =\,<\xi _{wh},\xi _{wh}>\, \\ 
& =\,1=\,<\xi _{h},\xi _{h}> \\ 
& =\,<\xi _{h},L_{w}^{*}\xi _{wh}>,
\end{array}
\,$
\end{center}

\strut

where $<,>$ is the inner product on the graph Hilbert space $H_{G}.$ Of
course, in the above formula we need the admissibility of $w$ and $h$ in $%
\mathbb{F}^{+}(G).$ However, even though $w$ and $h$ are not admissible
(i.e., $wh\notin \mathbb{F}^{+}(G)$), by the definition of $L_{w}^{\ast },$
we have that

\strut

\begin{center}
$
\begin{array}{ll}
<L_{w}\xi _{h},\xi _{h}> & =\,<0,\xi _{h}> \\ 
& =0=\,<\xi _{h},0> \\ 
& =\,<\xi _{h},L_{w}^{*}\xi _{h}>.
\end{array}
\,\,$
\end{center}

\strut

Notice that the creation operator $L$ and the annihilation operator $L^{*}$
satisfy that

\strut

(1.1) \ \ \ $L_{w}^{*}L_{w}=L_{y}$ \ \ and \ \ $L_{w}L_{w}^{*}=L_{x},$ \ for
all \ $w=xwy\in \mathbb{F}^{+}(G),$

\strut

\textbf{under the weak topology}, where $x,y\in V(G).$ Remark that if we
consider the von Neumann algebra $W^{*}(\{L_{w}\})$ generated by $L_{w}$ and 
$L_{w}^{*}$ in $B(H_{G}),$ then the projections $L_{y}$ and $L_{x}$ are
Murray-von Neumann equivalent, because there exists a partial isometry $%
L_{w} $ satisfying the relation (1.1). Indeed, if $w=xwy$ in $\mathbb{F}%
^{+}(G), $ with $x,y\in V(G),$ then under the weak topology we have that

\strut

(1,2) \ \ \ $L_{w}L_{w}^{*}L_{w}=L_{w}$ \ \ and \ \ $%
L_{w}^{*}L_{w}L_{w}^{*}=L_{w}^{*}.$

\strut

So, the creation operator $L_{w}$ is a partial isometry in $W^{*}(\{L_{w}\})$
in $B(H_{G}).$ Assume now that $v\in V(G).$ Then we can regard $v$ as $%
v=vvv. $ So,

\strut

(1.3) $\ \ \ \ \ \ \ \ \ L_{v}^{*}L_{v}=L_{v}=L_{v}L_{v}^{*}=L_{v}^{*}.$

\strut

This relation shows that $L_{v}$ is a projection in $B(H_{G})$ for all $v\in
V(G).$

\strut

Define the \textbf{graph }$W^{*}$\textbf{-algebra} $W^{*}(G)$ by

\strut

\begin{center}
$W^{*}(G)\overset{def}{=}\overline{%
\mathbb{C}[\{L_{w},L_{w}^{*}:w\in
\mathbb{F}^{+}(G)\}]}^{w}.$
\end{center}

\strut

Then all generators are either partial isometries or projections, by (1.2)
and (1.3). So, this graph $W^{*}$-algebra contains a rich structure, as a
von Neumann algebra. (This construction can be the generalization of that of
group von Neumann algebra.)\strut

\strut

\strut

\section{Lattice Path Model}

\strut

\strut

\strut

Throughout this section, let $G$ be a countable directed graph and let $%
W^{*}(G)$ be the graph $W^{*}$-algebra. Let $w_{1},$ $...,$ $w_{n}$ $\in $ $%
\Bbb{F}^{+}(G)$ and let $L_{w_{1}}^{u_{w_{1}}}$ $...$ $L_{w_{n}}^{u_{w_{n}}}$
$\in $ $W^{*}(G)$ be an operator in $W^{*}(G)$. In this section, we will
define a lattice path model for the random variable $%
L_{w_{1}}^{u_{w_{1}}}...L_{w_{n}}^{u_{w_{n}}}.$ Recall that if $w$ $=$ $%
e_{1} $ $...$ $e_{k}$ $\in $ $FP(G)$ with $e_{1},$ $...,$ $e_{k}$ $\in $ $%
E(G),$ then we can define the length $\left| w\right| $ of $w$ by $k.$ i.e.,
the length $\left| w\right| $ of $w$ is the cardinality $k$ of the
admissible edges $e_{1},...,e_{k}.$

\strut

\begin{definition}
Let $G$ be a countable directed graph and $\Bbb{F}^{+}(G),$ the free
semigroupoid. If $w\in \Bbb{F}^{+}(G),$ then $L_{w}$ is the corresponding $%
D_{G}$-valued random variable in $\left( W^{*}(G),E\right) .$ We define the
lattice path $l_{w}$ of $L_{w}$ and the lattice path $l_{w}^{-1}$ of $%
L_{w}^{*}$ by the lattice paths satisfying that ;

\strut

(i) \ \ the lattice path $l_{w}$ starts from $*=(0,0)$ on the $\Bbb{R}^{2}$%
-plane.

\strut

(ii) \ if $w\in V(G),$ then $l_{w}$ has its end point $(0,1).$

\strut

(iii) if $w\in E(G),$ then $l_{w}$ has its end point $(1,1).$

\strut

(iv) if $w\in E(G),$ then $l_{w}^{-1}$ has its end point $(-1,-1).$

\strut

(v) \ if $w\in FP(G)$ with $\left| w\right| =k,$ then $l_{w}$ has its end
point $(k,k).$

\strut

(vi) if $w\in FP(G)$ with $\left| w\right| =k,$ then $l_{w}^{-1}$ has its
end point $(-k,-k).$

\strut

Assume that finite paths $w_{1},...,w_{s}$ in $FP(G)$ satisfy that $%
w_{1}...w_{s}\in FP(G).$ Define the lattice path $l_{w_{1}...w_{s}}$ by the
connected lattice path of the lattice paths $l_{w_{1}},$ ..., $l_{w_{s}}.$
i.e.e, $l_{w_{2}}$ starts from $(k_{w_{1}},k_{w_{1}})\in \Bbb{R}^{+}$ and
ends at $(k_{w_{1}}+k_{w_{2}},k_{w_{1}}+k_{w_{2}}),$ where $\left|
w_{1}\right| =k_{w_{1}}$ and $\left| w_{2}\right| =k_{w_{2}}.$ Similarly, we
can define the lattice path $l_{w_{1}...w_{s}}^{-1}$ as the connected path
of $l_{w_{s}}^{-1},$ $l_{w_{s-1}}^{-1},$ ..., $l_{w_{1}}^{-1}.$
\end{definition}

\strut

\begin{definition}
Let $G$ be a countable directed graph and assume that $%
L_{w_{1}},...,L_{w_{n}}$ are generators of $W^{*}(G).$ Then we have the
lattice paths $l_{w_{1}},$ ..., $l_{w_{n}}$ of $L_{w_{1}},$ ..., $L_{w_{n}},$
respectively in $\Bbb{R}^{2}.$ Suppose that $%
L_{w_{1}}^{u_{w_{1}}}...L_{w_{n}}^{u_{w_{n}}}$ $\neq $ $0_{D_{G}}$ in $%
W^{*}(G),$ where $u_{w_{1}},...,u_{w_{n}}\in \{1,*\}.$ Define the lattice
path $l_{w_{1},...,w_{n}}^{u_{w_{1}},...,u_{w_{n}}}$ of nonzero $%
L_{w_{1}}^{u_{w_{1}}}...L_{w_{n}}^{u_{w_{n}}}$ by the connected lattice path
of $l_{w_{1}}^{t_{w_{1}}},$ ..., $l_{w_{n}}^{t_{w_{n}}},$ where $t_{w_{j}}=1$
if $u_{w_{j}}=1$ and $t_{w_{j}}=-1$ if $u_{w_{j}}=*.$ Assume that $%
L_{w_{1}}^{u_{w_{1}}}...L_{w_{n}}^{u_{w_{n}}}$ $=$ $0_{D_{G}}.$ Then the
empty set $\emptyset $ in $\Bbb{R}^{2}$ is the lattice path of it. We call
it the empty lattice path. By $LP_{n},$ we will denote the set of all
lattice paths of the $D_{G}$-valued random variables having their forms of $%
L_{w_{1}}^{u_{w_{1}}}...L_{w_{n}}^{u_{w_{n}}},$ including empty lattice path.
\end{definition}

\strut

Also, we will define the following important property on the set of all
lattice paths ;

\strut

\begin{definition}
Let $l_{w_{1},...,w_{n}}^{u_{w_{1}},...,u_{w_{n}}}\neq \emptyset $ be a
lattice path of $L_{w_{1}}^{u_{w_{1}}}...L_{w_{n}}^{u_{w_{n}}}\neq 0_{D_{G}}$
in $LP_{n}.$ If the lattice path $%
l_{w_{1},...,w_{n}}^{u_{w_{1}},...,u_{w_{n}}}$ starts from $*$ and ends on
the $*$-axis in $\Bbb{R}^{+},$ then we say that the lattice path $%
l_{w_{1},...,w_{n}}^{u_{w_{1}},...,u_{w_{n}}}$ has the $*$-axis-property. By 
$LP_{n}^{*},$ we will denote the set of all lattice paths having their forms
of $l_{w_{1},...,w_{n}}^{u_{w_{1}},...,u_{w_{n}}}$ which have the $*$%
-axis-property. By little abuse of notation, sometimes, we will say that the 
$D_{G}$-valued random variable $%
L_{w_{1}}^{u_{w_{1}}}...L_{w_{n}}^{u_{w_{n}}} $satisfies the $*$%
-axis-property if the lattice path $%
l_{w_{1},...,w_{n}}^{u_{w_{1}},...,u_{w_{n}}}$ of it has the $*$%
-axis-property.
\end{definition}

\strut

The following theorem shows that finding $E\left(
L_{w_{1}}^{u_{w_{1}}}...L_{w_{n}}^{u_{w_{n}}}\right) $ is checking the $*$%
-axis-property of $L_{w_{1}}^{u_{w_{1}}}...L_{w_{n}}^{u_{w_{n}}}.$

\strut

\begin{theorem}
Let $L_{w_{1}}^{u_{w_{1}}}...L_{w_{n}}^{u_{w_{n}}}\in W^{*}(G)$ be an
operator, where $u_{w_{1}},...,u_{w_{n}}\in \{1,*\}.$ Then $%
L_{w_{1}}^{u_{w_{1}}}...L_{w_{n}}^{u_{w_{n}}}$ $=$ $L_{v}$ if and only if $%
L_{w_{1}}^{u_{w_{1}}}...L_{w_{n}}^{u_{w_{n}}}$ has the $*$-axis-property.
\end{theorem}

\strut

\begin{proof}
($\Leftarrow $) Let $l=l_{w_{1},...,w_{n}}^{u_{w_{1}},...,u_{w_{n}}}\in
LP_{n}^{*}.$ Suppose that $w_{1}=vw_{1}v_{1}^{\prime }$ and $%
w_{n}=v_{n}w_{n}v_{n}^{\prime },$ for $v_{1},$ $v_{1}^{\prime },$ $v_{n},$ $%
v_{n}^{\prime }$ $\in $ $V(G).$ If $l$ is in $LP_{n}^{*},$ then

\strut

(2.1)$\ \ \ \ \ \ \ \left\{ 
\begin{array}{lll}
v_{1}=v_{n}^{\prime } &  & \text{if }u_{w_{1}}=1\text{ and }u_{w_{n}}=1 \\ 
&  &  \\ 
v_{1}=v_{n} &  & \text{if }u_{w_{1}}=1\text{ and }u_{w_{n}}=* \\ 
&  &  \\ 
v_{1}^{\prime }=v_{n}^{\prime } &  & \text{if }u_{w_{1}}=*\text{ and }%
u_{w_{n}}=1 \\ 
&  &  \\ 
v_{1}^{\prime }=v_{n} &  & \text{if }u_{w_{1}}=*\text{ and }u_{w_{n}}=*.
\end{array}
\right. $

\strut

By the definition of $LP_{n}^{*}$ and by (2.1), if $%
l_{w_{1},...,w_{n}}^{u_{w_{1}},...,u_{w_{n}}}\in LP_{n}^{*},$ then there
exists $v\in V(G)$ such that

\strut

$\ \ \ \ \ \ \ \ \ \ \ \ \ \ \ \ \ \ \ \ \ \ \ \
L_{w_{1}}^{u_{w_{1}}}...L_{w_{n}}^{u_{w_{n}}}=L_{v},$

where

(2.2) \ \ \ \ $\ \ \left\{ 
\begin{array}{ll}
v=v_{1}=v_{n}^{\prime } & \text{if }u_{w_{1}}=1\text{ and }u_{w_{n}}=1 \\ 
&  \\ 
v=v_{1}=v_{n} & \text{if }u_{w_{1}}=1\text{ and }u_{w_{n}}=* \\ 
&  \\ 
v=v_{1}^{\prime }=v_{n}^{\prime } & \text{if }u_{w_{1}}=*\text{ and }%
u_{w_{n}}=1 \\ 
&  \\ 
v=v_{1}^{\prime }=v_{n} & \text{if }u_{w_{1}}=*\text{ and }u_{w_{n}}=*.
\end{array}
\right. $

\strut

This shows that $L_{w_{1}}^{u_{w_{1}}}...L_{w_{n}}^{u_{w_{n}}}=L_{v}.$

\strut

($\Rightarrow $) Assume that $%
L_{w_{1}}^{u_{w_{1}}}...L_{w_{n}}^{u_{w_{n}}}=L_{v},$ in $W^{*}(G).$ Let $%
l=l_{w_{1},...,w_{n}}^{u_{w_{1}},...,u_{w_{n}}}\in LP_{n}$ be the lattice
path of the $D_{G}$-valued random variable $%
L_{w_{1}}^{u_{w_{1}}}...L_{w_{n}}^{u_{w_{n}}}.$ Trivially, $l\neq \emptyset
, $ since $l$ should be the connected lattice path. Assume that this
nonempty lattice path $l$ is contained in $LP_{n}\,\,\setminus \,LP_{n}^{*}.$
Then, under the same conditions of (2.1), we have that

\strut

(2.4) \ \ \ \ \ \ $\left\{ 
\begin{array}{lll}
v_{1}\neq v_{n}^{\prime } &  & \text{if }u_{w_{1}}=1\text{ and }u_{w_{n}}=1
\\ 
&  &  \\ 
v_{1}\neq v_{n} &  & \text{if }u_{w_{1}}=1\text{ and }u_{w_{n}}=* \\ 
&  &  \\ 
v_{1}^{\prime }\neq v_{n}^{\prime } &  & \text{if }u_{w_{1}}=*\text{ and }%
u_{w_{n}}=1 \\ 
&  &  \\ 
v_{1}^{\prime }\neq v_{n} &  & \text{if }u_{w_{1}}=*\text{ and }u_{w_{n}}=*.
\end{array}
\right. $

\strut

Therefore, by (2.2), there is no vertex $v$ satisfying $%
L_{w_{1}}^{u_{w_{1}}}...L_{w_{n}}^{u_{w_{n}}}=L_{v}.$ This contradict our
assumption.
\end{proof}

\strut

\strut

\strut

\section{Shadowed Graph Algebras}

\strut

\strut

\strut

In the previous chapter, we found one condition when the operator $%
T=L_{w_{1}}^{u_{w_{1}}}...L_{w_{n}}^{u_{w_{n}}}$ is a projection. In this
chapter, we will find the equivalent concept. To do that we need to define
the following new combinatorial object ;

\strut \strut

\begin{definition}
Let $G$ be a countable directed graph. Define the shadow $G^{-1}$ of the
graph $G$ by a graph with

\strut

$\ \ \ \ V(G^{-1})=V(G)$ \ and \ $E(G^{-1})=\{e^{-1}:e\in E(G)\}$,

\strut

where $e^{-1}$ is the reversely directed edge of $e.$ Then we can define the
shadowed graph $G^{\symbol{94}}$ $=$ $G$ $\cup $ $G^{-1}$ of $G$ by a graph
with

\strut

$\ \ \ \ \ \ \ V(G^{\symbol{94}})=V(G)$ \ and \ $E(G^{\symbol{94}})=E(G)\cup
E(G^{-1}).$

\strut

Let $\Bbb{F}^{+}(G^{\symbol{94}})$ be the free semigroupoid of the shadowed
graph $G^{\symbol{94}}.$ Define the relation $\mathbf{R}$ on $\Bbb{F}^{+}(G^{%
\symbol{94}})$ by

\strut

\ \ \ \ \ \ \ \ \ \ $\mathbf{R}$ $:$ \ $\ \ \ \ \ w^{-1}w=v^{\prime }$ \ \
and \ \ $ww^{-1}=v,$

\strut

whenever $w=vwv^{\prime }\in FP(G),$ with $v,v^{\prime }\in V(G).$ Define
the free groupoid $\Bbb{F}(G)$ of $G$ by the quotient set $\Bbb{F}^{+}(G^{%
\symbol{94}})\,/\,\mathbf{R}.$
\end{definition}

\bigskip

\begin{definition}
Let $G$ be a countable directed graph and let $G^{\symbol{94}}$ be the
shadowed graph of $G.$ Also, let $\Bbb{F}(G)$ be the free groupoid of $G.$
Let $l^{2}\left( \Bbb{F}(G)\right) $ be the Hilbert space with its Hilbert
basis $\{\xi _{w}:w\in \Bbb{F}(G)\}.$ Then we can define the creation
operator $M_{w}$ of $\xi _{w},$ as usual, under the admissibility with the
relation $\mathbf{R}.$ Define a weak-closed algebra $A\lg _{w}(G)$ by

\strut \strut

$\ \ \ \ \ \ \ \ \ \ \ A\lg _{w}(G)\overset{def}{=}\overline{\Bbb{C}%
[\{M_{w}:w\in \Bbb{F}(G)\}]}^{w}.$

\strut

This algebra $A\lg _{w}(G)$ is called the shadowed graph algebra of $G.$
\end{definition}

\strut

We will find the connection between our graph $W^{*}$-algebra $W^{*}(G)$ and
the shadowed graph algebra $A\lg _{w}(G)$. Notice that we can understand the
algebra $A\lg _{w}(G)$ as a weak-closed $*$-algebra with the involution

\strut

\begin{center}
$M_{w}^{*}=M_{w^{-1}},$ \ \ for all \ \ $w\in \Bbb{F}(G).$
\end{center}

\strut

This admit us to conclude that

\strut

\begin{theorem}
As weak-closed $*$-algebras, the graph $W^{*}$-algebra $W^{*}(G)$ and the
shadowed graph algebra $A\lg _{w}(G)$ are isomorphic. $\square $
\end{theorem}

\strut

The above theorem is easily proved because we can take the
generator-preserving linear map,

\strut

\begin{center}
$L_{w}^{u_{w}}\in W^{*}(G)\longmapsto M_{w^{t_{w}}}\in A\lg _{w}(G)$
\end{center}

\strut

where

\begin{center}
$t_{w}=\left\{ 
\begin{array}{lll}
1 &  & \text{if }u_{w}=1 \\ 
&  &  \\ 
-1 &  & \text{if }u_{w}=*,
\end{array}
\right. $
\end{center}

\strut

for all $w\in \Bbb{F}^{+}(G).$

\strut \bigskip \strut

Let $L_{w_{1}}^{u_{1}}...L_{w_{n}}^{u_{n}}\in W^{*}(G)$ be an operator,
where $u_{j}\in \{1,*\}.$ We can regard $%
L_{w_{1}}^{u_{1}}...L_{w_{n}}^{u_{n}}$ as $%
M_{w_{1}^{t_{1}}}...M_{w_{n}^{t_{n}}},$ by the previous theorem, where

\bigskip

\begin{center}
$t_{j}=\left\{ 
\begin{array}{cc}
1 & \text{if }u_{j}=1 \\ 
-1 & ~~\text{if }u_{j}=\ast ,
\end{array}
\right. $
\end{center}

\bigskip

for all \ $j=1,...,n,$ in $A\lg _{w}(G).$ Notice that, by definition,

\strut

\begin{center}
$M_{w_{1}^{t_{1}}}...M_{w_{n}^{t_{n}}}=M_{w_{1}^{t_{1}}...w_{n}^{t_{n}}}$
\end{center}

and

\begin{center}
$w_{1}^{t_{1}}...w_{n}^{t_{n}}\in \Bbb{F}(G)$ \ or \ $%
w_{1}^{t_{1}}...w_{n}^{t_{n}}\notin \Bbb{F}(G).$
\end{center}

\strut

By the previous discussion, we can have the following theorem ;

\strut

\begin{theorem}
Let $L_{w_{1}}^{u_{1}}...L_{w_{n}}^{u_{n}}\in W^{*}(G)$ be an operator,
where $u_{1},...,u_{n}\in \{1,*\}$, and let $%
M_{w_{1}^{t_{1}}...w_{n}^{t_{n}}}\in A\lg _{w}(G)$ be the corresponding
operator in $A\lg _{w}(G),$ where $t_{j}=1$ if $u_{j}=1$ and $t_{j}=-1$ if $%
u_{j}=*.$ Then $L_{w_{1}}^{u_{1}}...L_{w_{n}}^{u_{n}}=L_{v},$ for some $v\in
V(G),$ if and only if $w_{1}^{t_{1}}...w_{n}^{t_{n}}=v.$
\end{theorem}

\strut

\begin{proof}
It is clear from the definition of the free groupoid.
\end{proof}

\strut

\strut

\strut

\begin{quote}
\textbf{Reference}

\strut

\strut

{\small [1] I. Cho, Graph }$W^{*}${\small -Probability Spaces, (2004),
Preprint. }

{\small [2] D.W. Kribs and S.C. Power, Free Semigroupoid Algebras, (2003),
Preprint. }
\end{quote}

\end{document}